\documentstyle{amltd2004}
\begin{document}
\annalsline{155}{2002}
\received{December 19, 2000}
\startingpage{915}
\def\bye{\end{document}}
 \font\tenrm=cmr10

\def\eqref#1{(\ref{#1})}

\catcode`\@=11
\font\twelvemsb=msbm10 scaled 1100
\font\tenmsb=msbm10
\font\ninemsb=msbm10 scaled 800
\newfam\msbfam
\textfont\msbfam=\twelvemsb  \scriptfont\msbfam=\ninemsb
  \scriptscriptfont\msbfam=\ninemsb
\def\msb@{\hexnumber@\msbfam}
\def\Bbb{\relax\ifmmode\let\next\Bbb@\else
 \def\next{\errmessage{Use \string\Bbb\space only in math
mode}}\fi\next}
\def\Bbb@#1{{\Bbb@@{#1}}}
\def\Bbb@@#1{\fam\msbfam#1}
\catcode`\@=12

 \catcode`\@=11
\font\twelveeuf=eufm10 scaled 1100
\font\teneuf=eufm10
\font\nineeuf=eufm7 scaled 1100
\newfam\euffam
\textfont\euffam=\twelveeuf  \scriptfont\euffam=\teneuf
  \scriptscriptfont\euffam=\nineeuf
\def\euf@{\hexnumber@\euffam}
\def\frak{\relax\ifmmode\let\next\frak@\else
 \def\next{\errmessage{Use \string\frak\space only in math
mode}}\fi\next}
\def\frak@#1{{\frak@@{#1}}}
\def\frak@@#1{\fam\euffam#1}
\catcode`\@=12


\input amssym.def
\input amssym.tex

 \newcommand{\ad}{\mathop{\rm ad}\nolimits}
\newcommand{\Aut}{\mathop{\rm Aut}\nolimits}
\newcommand{\der}{\mathop{\rm ss}\nolimits}
\newcommand{\diag}{\mathop{\rm diag}\nolimits}
\newcommand{\disc}{\mathop{\rm disc}\nolimits}
\newcommand{\End}{\mathop{\rm End}\nolimits}
\newcommand{\Gal}{\mathop{\rm Gal}\nolimits}
\newcommand{\Lie}{\mathop{\rm Lie}\nolimits}
\newcommand{\Spin}{\mathop{\rm Spin}\nolimits}
\newcommand{\Tr}{\mathop{\rm Tr}\nolimits}

\newcommand{\bZ}{{\bf Z}}
\newcommand{\bQ}{{\bf Q}}
\newcommand{\bR}{{\bf R}}
\newcommand{\bC}{{\bf C}}
\newcommand{\bH}{{\bf H}}
\newcommand{\bS}{{\bf S}}

\newcommand{\sA}{{\cal A}}
\newcommand{\sB}{{\cal B}}
\newcommand{\sG}{{\cal G}}
\newcommand{\sV}{{\cal V}}
\newcommand{\sW}{{\cal W}}
\newcommand{\sY}{{\cal Y}}

\newcommand{\germg}{{\frak g}}
\newcommand{\germh}{{\frak h}}
\newcommand{\germk}{{\frak k}}
\newcommand{\germl}{{\frak l}}
\newcommand{\germo}{{\frak o}}
\newcommand{\germp}{{\frak p}}
\newcommand{\germs}{{\frak s}}
\newcommand{\germS}{{\frak S}}

\title{Hodge structures on abelian varieties\\ of type III}

\shorttitle{Abelian varieties of Type III} 

 \acknowledgements{2000 {\it Mathematics Subject Classification}.  Primary 14C30, 14K20; Secondary 20G30.}
\author{Salman Abdulali}\institutions{East Carolina University,
Greenville, NC\\
{\eightpoint {\it E-mail address\/}: abdulalis@mail.ecu.edu}}
 
 \centerline{\bf Abstract}
\vglue4pt
We show that the usual Hodge conjecture implies the general Hodge
conjecture for certain abelian varieties of type III, and use this to
deduce the general Hodge conjecture for all powers of certain
$4$-dimensional abelian varieties of type III. We also show the
existence of a Hodge structure $M$ such that $M$ occurs in the
cohomology of an abelian variety, but the Tate twist $M(1)$ does not
occur in the cohomology of any abelian variety, even though it is
effective.
 \vglue-16pt
 \section{Introduction}
 \vglue-4pt
Consider a Hodge structure $V_{\bC} = \mathbold{\oplus}_{p+q=n} V^{p,q}$ of weight $n$.
The Hodge structure is said to be {\it effective} if $V^{p,q} = 0$
unless $p,q \geq 0$.
It is said to be {\it geometric} if it is isomorphic to a Hodge
substructure of the cohomology of a smooth projective algebraic
variety over $\bC$.
For an integer $m$, the {\it Tate twist} $V(m)$ is the Hodge
structure of weight $n-2m$ defined by $V(m)^{p,q} := V^{p+m,q+m}$.

The general Hodge conjecture as formulated by Grothendieck
\cite{Grothendieck} implies that any Tate twist of a geometric Hodge
structure is again geometric, provided that it is effective.
We say that a smooth projective variety $X$ is {\it dominated} by a
class of varieties $\sY$ if, for every irreducible Hodge structure
$M$ occurring in the cohomology of~$X$, the Tate twist $M(h)$ occurs
in the cohomology of some member of~$\sY$, where
$$h := \min \{p \mid M^{p,n-p} \neq 0 \}$$
is the {\it height} of the Hodge structure.
Grothendieck observed that the general Hodge conjecture for $X$ is
implied by the existence of a class of varieties $\sY$ which
dominates~$X$, together with the usual Hodge conjecture for $Y \times
X$ for each $Y \in \sY$ (cf.~\cite[Proposition 2.1,
p.~343]{Abdulali2}).

In a series of papers (\cite{Abdulali2}, \cite{Abdulali4},
\cite{Abdulali5}) we showed that certain abelian varieties are
dominated by subclasses of the class of all abelian varieties, and
used this to deduce the general Hodge conjecture in some cases.
This raised the question of whether every abelian variety is
dominated by the class of all abelian varieties.
We answer this question in the negative in this paper (Theorem
\ref{maintheorem}).

Let $A$ be a simple abelian variety of type III.
Recall that $A$ is said to be of type III if its endomorphism algebra
$D := \End(A) \otimes \bQ$ is a totally definite quaternion algebra
over a totally real number field.
Assume that $A$ is a general member of a PEL-family. In
\cite[Theorem~5.1, p.~348]{Abdulali2} we showed that $A$ is dominated
by the set of powers of itself, provided that $m := \dim_D H_1(A,
\bQ)$ is odd.
Using this we were able to prove the general Hodge conjecture for all
powers of certain $6$-dimensional abelian varieties of type III in
\cite[Example 5.1, p.~674]{Abdulali3}.
In this paper we show (Theorem \ref{nonsquare} and Proposition
\ref{squareprop}) that $A$ is dominated by the powers of itself if
and only if $\disc T$, the discriminant of the $D$-valued
skewhermitian form determined by a polarization of $A$, is not a
square.
As a consequence we prove the general Hodge conjecture for all powers
of certain $4$-dimensional abelian varieties of type III (Corollary
\ref{fourdim}).

When $\disc T$ is a square and $m=4$, we show that there exists an
abelian variety $B$ such that $A$ is dominated by abelian varieties
of the form $A^i \times B^j$ (Theorem \ref{d4}).
However, when $\disc T$ is a square and $m>4$, we show that there
exists a Hodge structure $M \subset H^m(A, \bQ)$ such that $M(1)$ is
effective but does not occur in the cohomology of any abelian variety
(Theorem \ref{maintheorem}).
The least dimension of the abelian varieties in which these
extraordinary Hodge structures occur is $12$, in which case $M
\subset H^6 (A, \bQ)$.
Weil \cite{Weil} proposed that certain Hodge cycles be considered as
possible counterexamples to the usual Hodge conjecture.
In the same spirit, we suggest that these extraordinary Hodge
structures be used as a testing ground for the general Hodge
conjecture.

\vglue4pt {\it Acknowledgements}.
Lemma~\ref{quotient} was formulated during a discussion with Chad
Schoen at Banff in June 1998. I thank him for this and other
suggestions. I thank the referees for their careful reading of the
manuscript; their comments have helped improve the paper.
\vglue4pt

{\it Notation and conventions}.
All abelian varieties are over $\bC$. For an abelian variety $A$, we
let $D(A) := \End(A) \otimes \bQ$ be the endomorphism algebra of $A$,
we denote by $G(A)$ the Hodge group of $A$ (\S 2.1), and we write
$\germg (A)$ for the Lie algebra of $G(A)$. The semisimple part
(derived group) of $G(A)$ is denoted by $G^{\der}(A)$, and
$\germg^{\der} (A) := \Lie G^{\der}(A)$ is its Lie algebra.

The definition of {\it dominated} given in this paper is slightly
different from that given in \cite[p.~344]{Abdulali2}; the two
definitions are equivalent if the general Hodge conjecture is true.
 
\vglue-6pt
\section{Abelian varieties and Hodge groups}
\vglue-6pt

 2.1. {\it Hodge groups}.
By a Hodge structure of weight $n$, we mean a finite-dimensional
vector space $V$ over $\bQ$, and a direct sum decomposition $V_{\bC} =
\mathbold{\oplus}_{p+q=n} V^{p,q}$  such that $\overline{V^{p,q}} = V^{q,p}$.
Let $\bS$ be the real algebraic group such that $\bS(\bR)$ is the
unit circle in the complex plane. Then a Hodge structure $V$
determines a morphism
$$\tau_V : \bS \to {\rm GL}(V_{\bR})$$
of real algebraic groups such that $\tau_V(e^{i\theta})$ acts as
multiplication by $e^{(p-q)i\theta}$ on $V^{p,q}$. The {\it Hodge
group} of $V$ is defined to be the smallest $\bQ$-subgroup, $G$, of
${\rm GL}(V)$ such that $G_{\bR}$ contains the image of $\tau_V$.

The {\it Hodge group}, $G(A)$, of an abelian variety $A$ is defined
to be the Hodge group of $V := H^1(A, \bQ)$.
A Riemann form for $A$ gives an isomorphism of $V$ with $V^\star$
which is a morphism of Hodge structures, so the Hodge group of $A$
may also be defined as the Hodge group of $H_1(A, \bQ) \cong H^1(A,
\bQ)^{\star}$.
The Hodge group of an abelian variety is a connected, reductive group
of hermitian type which has no nontrivial, connected, normal
$\bQ$-subgroup $H$ such that $H(\bR)$ is compact (\cite{Mumford1},
\cite{900}, \cite{Gordon}).

\specialnumber{2.1.1}\proclaim{Lemma}
\label{quotient}
Let $M$ be a Hodge structure contained in the cohomology of an
abelian variety $A${\rm .}
Then the Hodge group of $M$ is a quotient of the Hodge group of $A${\rm .}
\endproclaim

\demo{Proof}
Let $G$ be the Hodge group of $A$, $H$ the Hodge group of $M$, and,
$V := H^1(A, \bQ)$.
Since $M$ is a $G$-submodule of the exterior algebra of $V$, we have
a morphism $\varphi : G \to {\rm GL}(M)$. This is compatible with the Hodge
structures in the sense that
$\varphi (\tau_V (g)) = \tau_M (g)$ for all $g \in \bS(\bR)$.
Therefore $\varphi (G)$ contains the image of $\tau_M$, and hence
contains $H$. Since $\varphi^{-1}(H)$ is a $\bQ$-subgroup of $G$
containing the image of $\tau_V$, we see that $\varphi^{-1}(H) = G$.
Thus $H = \varphi(G)$.
\enddemo

2.2. {\it Kuga fiber varieties}.
A general reference for this subsection is Satake's book \cite{Satakebook}.

Let $G$ be a connected, semisimple, linear algebraic group over $\bQ$.
Assume that $G$ is of hermitian type, and has no nontrivial,
connected, normal $\bQ$-subgroup $H$ such that $H(\bR)$ is compact.
Then $X := G(\bR)^0/K$ is a bounded symmetric domain, where $K$ is a
maximal compact subgroup of $G(\bR)^0$.
Let $\germg := \Lie G$ be the Lie algebra of $G$, $\germk := \Lie K$,
and let $\germg_{\bR} = \germk \oplus \germp$ be the corresponding
Cartan decomposition.
Let $o$ be the unique fixed point of $K$ in $X$.
Differentiating the natural map $G(\bR)^0 \to X$ gives an isomorphism
of $\germp$ with $T_o(X)$, the tangent space of $X$ at $o$, and there
exists a unique $H_0 \in Z(\germk)$, called the  $H$-{\it element} at
$o$, such that $\ad H_0|\germp$ is the complex structure on $T_o(X)$.

Let $\beta$ be a nondegenerate alternating form on a
finite-dimensional rational vector space $V$. The symplectic group
${\rm Sp}(V,\beta)$ is a $\bQ$-algebraic group of hermitian type; the
associated symmetric domain is the {\it Siegel space}
\begin{eqnarray*}
\germS(V,\beta)& :=& \{J \in {\rm GL}(V_{\bR}) \mid J^2 = -I \hbox{ and } \\
&&\qquad\beta(x,Jy) \hbox{ is symmetric, positive definite}\}.
\end{eqnarray*}
${\rm Sp}(V,\beta)$ acts on $\germS(V,\beta)$ by conjugation.
The $H$-element at a point $J \in \germS(V,\beta)$ is $J/2$. 

Let $\rho : G \to {\rm Sp}(V, \beta)$ be a representation defined over
$\bQ$. We say that $\rho$ satisfies the {\it $H_1$-condition}
relative to the $H$-elements $H_0$ and $H'_0 = J/2$ if
$$[d \rho (H_0) - H'_0, d \rho (g)] = 0 \qquad \hbox{for all }g \in
\germg_{\bR},$$
and we say that $\rho$ satisfies the {\it $H_2$-condition} if
$$d \rho (H_0) = H'_0.$$

Assume that the $H_1$-condition is satisfied relative to $H_0$ and $H_0'$.
Then there exists a unique holomorphic map $\tau : X \to
\germS(V,\beta)$ such that $\tau (o) = J$, and the pair $(\rho,
\tau)$ is equivariant in the sense that
$$
\tau (g \cdot x) = \rho(g) \cdot \tau (x) \qquad \hbox{for all } g
\in G(\bR)^0, x \in X.
$$
Let $\Gamma$ be a torsion-free arithmetic subgroup of $G(\bQ)$.
Let $L$ be a lattice in $V$ such that $\rho(\Gamma)L = L$.
The quotient
$$\sA := (\Gamma \ltimes_\rho L) \backslash (X \times V_{\bR})$$
is a smooth quasiprojective algebraic variety (see
\cite[p.~74]{900}); the natural projection $\sA \to \sV := \Gamma
\backslash X$ makes $\sA$ a fiber variety over $\sV$ called a
{\it Kuga fiber variety}.
The fiber $\sA_P$ over any point $P \in \sV$ is an abelian variety
isomorphic to the torus $V_{\bR}/L$ with the complex structure $\tau
(x)$, where $x$ is a point in $X$ lying over $P$.

\specialnumber{2.2.1}\proclaim{Lemma}
\label{rigidity}
Let $\rho_j : G \to {\rm Sp}(V_j,\beta_j)${\rm ,} $j = 1, 2${\rm ,} be
$H_2$\/{\rm -}\/representations defining Kuga fiber varieties $\sA_j \to \sV${\rm .}
Let $P \in \sV${\rm .}
If $\rho_1$ and $\rho_2$ are equivalent as representations of $G${\rm ,}
then{\rm ,} the abelian varieties $\sA_{1,P}$ and $\sA_{2,P}$ are isogenous{\rm .}
\endproclaim

\demo{Proof}
Let $\psi : V_1 \to V_2$ be an equivalence, and $\psi_{\star}$ the
induced isomorphism from ${\rm GL}(V_1)$ to ${\rm GL}(V_2)$.
The abelian variety $\sA_{j,P}$ is the torus $V_{j,\bR}/L_j$ with the
complex structure $J_j$, where $L_j$ is a lattice in $V_j$, and $J_j
= 2 d\rho_j(H_0)$.
Since $\psi_{\star} (J_1) = J_2$, $\psi$ induces an isogeny from
$\sA_{1,P}$ to $\sA_{2,P}$.
\enddemo

2.3. {\it Hodge families}.
Given an abelian variety $A$, Mumford \cite{Mumford1} constructed a
Kuga fiber variety having $A$ as a general fiber.
These Kuga fiber varieties are called {\it Hodge families}, and have
interpretations as solutions to fine moduli problems.
We briefly review Mumford's construction here.

Let $A$ be an abelian variety, and $V := H_1(A, \bQ)$.
As a complex torus, $A$ is isomorphic to $V_{\bR} / V_{\bZ}$ with a
complex structure $J$, where $V_{\bZ}$ is the lattice $H_1(A, \bZ)$.
A polarization of $A$ determines a Riemann form, i.e., an alternating
form, $\beta$, on $V$ such that $\beta(V_{\bZ}, V_{\bZ}) \subset
\bZ$, and $\beta (x,Jy)$ is symmetric and positive definite.
Let $X$ be the symmetric domain belonging to $G^{\der}(A)$, and
choose a base point $o \in X$.
The inclusion
$$\rho : G^{\der}(A) \to {\rm Sp}(V, \beta)$$
satisfies the $H_1$-condition with respect to the $H$-elements at $o$ and $J$.
Let $\Gamma$ be a torsion-free arithmetic subgroup of $G^{\der}(A)$
such that $\rho(\Gamma)(V_{\bZ}) \subset V_{\bZ}$.
The Kuga fiber variety $\sA \to \sV$ obtained from this data is
called the {\it Hodge family determined by} $A$.
If $P$ is the image of $o$ in $\sV$, then the fiber $\sA_P$ is
isomorphic to $A$, and the Hodge group of every member of the family
is contained in $G(A)$.

\proclaimtitle{[1, Prop.~2.2, p.~1124]}
\specialnumber{2.3.1}\proclaim{Proposition}
\label{h2}
Let $G$ be a connected{\rm ,} semisimple{\rm ,} linear algebraic group over $\bQ$
of hermitian type{\rm .} Assume that $G$ has no nontrivial{\rm ,} connected{\rm ,}
normal $\bQ$\/{\rm -}\/subgroup $H$ such that $H(\bR)$ is compact{\rm .} Any Kuga
fiber variety $\sA \to \sV$ constructed from a symplectic
representation $\rho : G \to {\rm Sp}(V, \beta)$ satisfying the
$H_2$\/{\rm -}\/condition{\rm ,} is a Hodge family{\rm ,} and the Hodge group of a general
fiber equals $\rho(G)${\rm .}
\endproclaim

2.4. {\it Primary representations}.
Recall that a fully reducible representation is called {\it primary}
if all irreducible representations contained in it are equivalent.
Let $\rho : G \to {\rm Sp}(V, \beta)$ be an $H_1$-representation of a
semisimple hermitian group $G$ over $\bQ$. We have the primary
decomposition $V = \mathbold{\oplus}_{\alpha} V_{\alpha}$, where the
$V_{\alpha}$ are the maximal primary $G$-submodules of $V$
(cf.~Satake \cite[p.~167]{Satakebook}). Each representation
$$G \to {\rm Sp}(V_{\alpha}, \beta | V_{\alpha} \times V_{\alpha})$$
satisfies the $H_1$-condition (see Satake \cite[Lemma IV.4.2,
p.~181]{Satakebook} and the remarks following it).

\specialnumber{2.4.1}\proclaim{Lemma}
\label{primarylemma}
Let $\rho : G \to {\rm Sp}(V, \beta)$ be a primary representation
satisfying the $H_1$\/{\rm -}\/condition{\rm ,} and defined over $\bQ${\rm .} Suppose $G =
G_1 \times G_2${\rm ,} and $G_1(\bR)$ has no compact factors{\rm .} Then either
$G_1$ or $G_2$ acts trivially on $V${\rm .}
\endproclaim

\demo{Proof}
Suppose both $G_1$ and $G_2$ act nontrivially on $V$. Let $U$ be an
irreducible submodule of $V_{\overline{\bQ}}$. Then $U^{\tau}$ is
also an irreducible submodule of $V_{\overline{\bQ}}$ for any
automorphism $\tau$ of $\overline{\bQ}$, and we have
$V_{\overline{\bQ}}$ equal, up to multiplicity, to the sum of the
distinct $\Aut(\overline{\bQ})$-conjugates of $U$ (see Satake
\cite[\S 1.1, pp.~218--221]{Satake2}). If $G_1$ (resp.~$G_2$) acts
trivially on $U$, then $G_1$ (resp.~$G_2$) acts trivially on every
conjugate of $U$, so we infer that both $G_1$ and $G_2$ act
nontrivially on $U$.

Some simple factor $H$ of $G_2$ acts nontrivially on $U$. We have $H
= R_{K/\bQ}\widehat{H}$, where $K$ is a totally real number field,
$\widehat H$ is an absolutely simple group over $K$, and, $R_{K/\bQ}$
means restriction of scalars from $K$ to $\bQ$ (see Satake
\cite[p.~259]{ Satake2}). Let $S$ be the set of embeddings of $K$
into $\bR$, let $L$ be the Galois closure of $K$, and let
$\sG := \Gal(L/\bQ)$.
Then $\sG$ acts transitively on $S$. We have
$H_L = \prod_{a \in S} H_a$, where $H_a := \widehat{H} \otimes_{K,a} L$. Let
$$S_0 := \{a \in S \mid H_a(\bR) \hbox{ is not compact} \}.$$
One of the factors, say $H_a$, acts nontrivially on $U$. Since $\sG$
acts transitively on $S$, there exists $\sigma \in \sG$ such that
$a^{\sigma} = b$, and, $b \in S_0$. Extend $\sigma$ to an
automorphism of $\overline{\bQ}$, and denote the extension again by
$\sigma$. Then $U^{\sigma}$ is an irreducible subrepresentation of
$V_{\overline{\bQ}}$ such that both $G_1$ and $H_b$ act nontrivially
on it.

Let $V_{\bR} = \mathbold{\oplus}_{\alpha} V_{\alpha}$ be the primary
decomposition of $V_{\bR}$ as a $G_{\bR}$-module.
Then we have $V_{\bC} = \mathbold{\oplus}_{\alpha} V_{\alpha,\bC}$. Since
$U^\sigma_{\bC}$ is an irreducible submodule of $V_{\bC}$, there exists
$\alpha$ such that $V_{\alpha,\bC}$ contains a submodule equivalent
to $U^{\sigma}_{\bC}$. Then both $G_1$ and $H_b$ act nontrivially on
$V_\alpha$.
Since $G_1(\bR)$ has no compact factors, and  $H_b(\bR)$ is not
compact, this contradicts Satake's result that each $\bR$-primary
subrepresentation of $\rho_{\bR}$ can be nontrivial on only one
noncompact factor of $G(\bR)$ (Satake
\cite[pp.~185--186]{Satakebook}).
\enddemo

\section{An algebraic group and its Lie algebra}
\label{groupalgebra}

3.1. {\it Quaternionic skewhermitian forms}.
Let $k$ be a totally real number field, $D$ a totally definite
quaternion algebra over $k$, and $x \mapsto \bar{x}$ the canonical
involution of $D$.
Let $V$ be a left vector space of dimension $m$ over $D$. Let $T$ be
a $D$-valued skewhermitian form on $V$, i.e., a $k$-bilinear form $T
: V \times V \to D$ such that $T(ax,by) = aT(x, y)\bar{b}$ and $T(y,
x) = -\overline{T(x, y)}$ for all $a,b \in D$ and all $x,y \in V$.
We assume that $T$ is nondegenerate.
Recall that the {\it determinant} $\det T$ is the reduced norm from
$M_n(D)$ to $k$ of a matrix representation of $T$, and the
{\it discriminant} of $T$ is $\disc T = (-1)^m \det T$.
The determinant and discriminant are only well-defined modulo
squares; i.e., they are elements of $k^{\times} / k^{\times 2}$.

\demo{Example {\rm 3.1.1}}
Recall that an element $\lambda \in D$ is called a {\it pure}
quaternion if $\bar{\lambda} = -\lambda$.
The reduced norm of such an element is given by $N(\lambda) = -\lambda^2$.
Let $\lambda_1, \dots, \lambda_m$ be nonzero pure quaternions in $D$.
Then $\diag(\lambda_1, \dots, \lambda_m)$ is the matrix of a
skewhermitian form $T$ on $V = D^m$.
We have $\disc T = \prod_i \lambda_i^2$.
\enddemo

{\it Example} 3.1.2.
As a special case of the previous example, suppose $\lambda =
\lambda_1 = \cdots = \lambda_m$.
Then $\disc T = \lambda^{2m}$ is a square in $k$ if and only if $m$ is even.
 
\demo{Example {\rm 3.1.3}}
As another special case of Example~3.1.1, let $D$ be the
quaternion algebra over $\bQ$ with basis $\{ 1,i,j,k \}$, where $i^2
= -1$, $j^2 = -3$, and, $ij = k$.
Take $m \geq 2$, let $\lambda_1 = \lambda_2 = \cdots = \lambda_{m-1} =
i$, and let $\lambda_m = j$.
Then $\disc T = \pm 3$ is not a square.
\enddemo

3.2. {\it An algebraic group}.
Let $\widehat{G}$ be the connected component of the identity in the
group $\Aut_D(V, T)$. Then $\widehat{G}$ is an algebraic group over
$k$ of type $D_m$. It is semisimple for $m \geq 2$, and absolutely
simple for $m > 2$.

Let $G$ be the group $R_{k/\bQ}\widehat{G}$ obtained by restricting
scalars from $k$ to $\bQ$. It is a semisimple algebraic group over
$\bQ$ for $m \geq 2$, and is $\bQ$-simple for $m > 2$. Let $S$ be the
set of embeddings of $k$ into $\bR$. Then the real Lie group $G(\bR)$
is given by
\begin{equation}
\label{gproduct}
G(\bR) \cong \prod_{\alpha \in S} G_{\alpha} \speqnu{3.2.1}
\end{equation}
where $G_\alpha$ is the group of real points of
$\widehat{G} \otimes_{k,\alpha} \bR$. 
Each $G_{\alpha}$ is isomorphic to the group of type $D_m$ called
${\rm SO}^{\star}(2m)$ by Helgason \cite[p.~445]{Helgason} and
${\rm SU}^-(m,\bH)$ by Satake (see \cite[Exercise 2, p.~278]{Satakebook}).
The real rank of $G_{\alpha}$ is $[m/2]$. It has Tits index
${}^2D_{m,(m-1)/2}^{(2)}$ if $m$ is odd, and Tits index
${}^1D_{m,m/2}^{(2)}$ if $m$ is even \cite[pp.~56--57]{Tits}.

\demo{{\rm 3.3.} The Lie algebra}
Let $\alpha$ be an embedding of $k$ into $\bR$.
We now recall the description of the Lie algebra $\germg_{\alpha}$ of
$G_{\alpha}$, and its representations, following Satake \cite[\S 3.3,
pp.~449--451]{Satake1} and Varadarajan \cite[p.~329]{Varadarajan}.
We have
$$
\germg_\alpha = \left\{
\left( \begin{array}{cc} X_1 & X_{12}\\ -\overline{X}_{12} & \overline{X}_1
\end{array}\right)
\in \germs\germl_{2m}(\bC) \Biggm|
\begin{array}{c} X_1, X_{12} \in M_m(\bC) \\ {}^t\overline{X}_1 = -X_1,
{}^tX_{12} = -X_{12} \end{array}
\right\},
$$
a Cartan subalgebra is given by
$$\germh = \{ \diag (a_1, \dots , a_m, -a_1, \dots, -a_m) \in
\germs\germl_{2m}(\bC) \mid a_j \in i\bR \},$$
and an $H$-element is given by
$$H_0 = \left( \begin{array}{cc} \frac{i}{2} I_m & 0 \\ 0 & \frac{-i}{2} I_m
\end{array}\right).$$

Define $\lambda_j : \germh_{\bC} \to \bC$ by $\lambda_j (\diag (a_1,
\dots, a_{2m})) = a_j$.
A simple system of roots of $\germg_{\alpha, \bC}$ is given by
$\Delta = \{ \alpha_1, \dots, \alpha_m \}$, where $\alpha_j =
\lambda_j - \lambda_{j+1}$ for $1 \leq j \leq m-1$, and, $\alpha_m =
\lambda_{m-1} + \lambda_m$.

The fundamental weights of $\germg_{\alpha, \bC}$ with respect to
$\germh_{\bC}$ are
$\mu_1, \dots , \mu_m$, where
\begin{eqnarray*}
\mu_j &=& \lambda_1 + \cdots + \lambda_j \qquad \hbox{for } 1 \leq j \leq m-2,\\
\mu_{m-1} &=& (\lambda_1 + \cdots + \lambda_{m-1} - \lambda_m)/2, \\
\mu_m &=& (\lambda_1 + \cdots + \lambda_{m-1} + \lambda_m)/2.
\end{eqnarray*}

Let $W$ be the standard representation of $\germg_{\alpha, \bC}$.
Then $\bigwedge^j W$ is an irreducible $\germg_{\alpha, \bC}$-module
with highest weight $\mu_j$ for $1 \leq j \leq m-2$,
$\bigwedge^{m-1}W$ is an irreducible $\germg_{\alpha, \bC}$-module
with highest weight $\mu_{m-1} + \mu_m$, while $\bigwedge^m W$ splits
as the direct sum of two irreducible $\germg_{\alpha,
\bC}$-submodules, $U_1$ and $U_2$, having highest weights
$2\mu_{m-1}$ and $2\mu_m$, respectively \cite[Chap.~4, Exercise 23,
p.~394]{Varadarajan}.

\specialnumber{3.3.1}\proclaim{Lemma}
\label{notappear}
Let ${\rm SO}_{2m}(\bC) \to {\rm GL}(W)$ be the standard representation of
${\rm SO}_{2m}(\bC)${\rm ,} where $m \geq 2${\rm .} Then{\rm ,} the irreducible
representations with highest weights $2\mu_{m-1}$ and $2\mu_m$ do not
occur in $\bigwedge^{m-2} (W^n)${\rm ,} for any positive integer~$n${\rm .}
\endproclaim

\demo{Proof}
An irreducible representation which occurs in $\bigwedge^{m-2} (W^n)$
must occur in $W^{\mathbold{\scriptstyle\otimes} {m-2}}$.
The set of weights of the standard representation of
$\germs\germo_{2m}(\bC)$ on $W$ is
$$\Lambda = \{ \lambda_1, \dots, \lambda_m, -\lambda_1, \dots, -\lambda_m \}.$$
The weights of the representation $W^{\otimes {m-2}}$ are of the form
$a_1 + \cdots + a_{m-2}$, where $a_j \in \Lambda$.
Therefore
$$2\mu_{m-1} = \lambda_1 + \cdots + \lambda_{m-1} - \lambda_m$$
and
$$2\mu_m = \lambda_1 + \cdots + \lambda_{m-1} + \lambda_m$$
are not weights of $W^{ \mathbold{\scriptstyle\otimes} {m-2}}$.
\enddemo

\section{Abelian varieties of type III}

Let $A$ be a simple abelian variety, and $D := D(A) := \End(A)
\otimes \bQ$ its endomorphism algebra.
We assume in this section that $A$ is of type~III; i.e., $D$ is a
totally definite quaternion algebra over a totally real field~$k$.
Let $V := H_1(A, \bQ)$.
Then $V$ is a left vector space over $D$; let $m$ be its dimension.
By Shimura \cite[Proposition~15, p.~177]{Shimura1} we have $m \geq 2$.

Let $\beta : V \times V \to \bQ$ be a Riemann form for $A$.
There exists a unique $D$-valued skewhermitian form $T$ on $V$ such
that $\beta (x,y) = \Tr_{D/\bQ} T(x,y)$ for all $x,y \in V$, where
$\Tr_{D/\bQ}$ is the reduced trace from $D$ to $\bQ$ (Shimura
\cite[Lemma~1.2, p.~162]{Shimura2}).
The group $\widehat{G} := \Aut_{D}(V,T)^0$ is then of the type
discussed in \S 3.2, and we shall henceforth use the
notations introduced in \S \ref{groupalgebra}. The group $G$ is a
semisimple group of hermitian type, and its inclusion in ${\rm Sp}(V,
\beta)$ satisfies the $H_2$-condition.
The resulting Kuga fiber variety is a PEL-family in the sense of
Shimura \cite{Shimura3}.
For proofs of these facts, and more details, we refer the reader to
Satake \cite[Example~IV.7.2, p.~200]{Satakebook}, and to
\cite{Abdulali3}.

\specialnumber{4.1} \proclaim{Theorem}
\label{nonsquare}
Let $A$ be a simple abelian variety of type {\rm III.}
Assume that $A$ is a general member of a {\rm PEL-}\/family{\rm .}
Let $T$ be the skewhermitian form determined by a polarization of $A${\rm .}
If the discriminant of $T$ is not a square{\rm ,} then any power of $A$ is
dominated by the set of powers of $A${\rm ,} and the usual Hodge conjecture
for $A$ implies the general Hodge conjecture for all powers of~$A${\rm .}
\endproclaim

\demo{Proof}
Proposition \ref{h2} implies that the Hodge group of $A$ is $G$.
We shall now explicitly describe the Hodge structure on the
cohomology of $A$, in terms of the action of $G$ on $V$ (see also the
proof of \cite[Theorem 5.1, pp.~348--353]{Abdulali2}).
Corresponding to the decomposition \eqref{gproduct} of $G(\bR)$, we have
$V_{\bR} \cong \mathbold{\oplus}_{\alpha \in S} V_{\alpha}$,
where each $V_{\alpha}$ is a real Hodge substructure of $V_{\bR}$
such that $G_{\beta}$ acts trivially on $V_{\alpha}$ unless $\alpha =
\beta$.
The action of $\germg_{\alpha}$ on $V_{\alpha}$ is given by
$$
\left( \begin{array}{cc} X_1 & X_{12}\\ -\overline{X}_{12} & \overline{X}_1
\end{array}\right)
\mapsto
\left( \begin{array}{cccc}
X_1 & 0 & 0 & -X_{12}\\
0 & X_1 & X_{12} & 0\\
0 & -\overline{X}_{12} & \overline{X}_1 & 0\\
\overline{X}_{12} & 0 & 0 & \overline{X}_1
\end{array}\right)
$$
(Satake \cite[(17), p.~433 and \S 3.3, pp.~449--451]{Satake1}) with
respect to a suitable basis
of $V_{\alpha,\bC}$, which we call
$$\sB = \{v_{m+1}, \dots, v_{2m}, u_1, \dots, u_{2m}, v_1, \dots, v_m\}.$$
This notation is somewhat perverse, but has been chosen to simplify
some of the expressions later in the proof.
We observe that
$$
V_{\alpha,\bC} = W_{\alpha} \oplus \overline{W}_{\alpha}, \speqnu{4.1.2}
$$
where $W_{\alpha}$ has basis $\{u_1, \dots, u_{2m} \}$, and its
complex conjugate, $\overline{W}_\alpha$, has basis $\{v_{m+1},
\dots, v_{2m}, v_1, \dots, v_m \}$.
The action of $\germg_{\alpha,\bC} \cong \germs\germo_{2m}(\bC)$ on
both $W_{\alpha}$ and $\overline{W}_{\alpha}$ is equivalent to the
standard representation.

With respect to the basis $\sB$, the complex structure on
$V_{\alpha}$ is given by
$$J = 2H_0 = \diag (iI_{2m}, -iI_{2m})$$
(Satake \cite[p.~431]{Satake1}).
  From this we see that the Hodge structure of $V_{\alpha}$ is
$$V_{\alpha,\bC} = V^{1,0} \oplus V^{0,1},$$
where $V^{1,0}$ has the basis $\{v_{m+1}, \dots, v_{2m}, u_1, \dots, u_m \}$,
and, $V^{0,1}$ has the basis $\{u_{m+1}, \dots, u_{2m}, v_1, \dots, v_m\}$.
It follows that
$$W_{\alpha} = W_{\alpha}^{1,0} \oplus W_{\alpha}^{0,1},$$
where $W_{\alpha}^{1,0}$ has basis $\{u_1, \dots, u_m\}$ and consists
of $(1,0)$-forms, while $W_{\alpha}^{0,1}$ has basis $\{u_{m+1},
\dots, u_{2m}\}$ and consists of $(0,1)$-forms.

For $1 \leq j \leq m-2$, $\bigwedge^j W_{\alpha}$ is an irreducible
$G_{\alpha} (\bC)$-module with highest weight $\mu_j$.
It contains the $(j,0)$-form $u_1 \wedge \cdots \wedge u_j$.
The irreducible $G_{\alpha} (\bC)$-module $\bigwedge^{m-1}W_{\alpha}$
has highest weight $\mu_{m-1} + \mu_m$.
It contains the $(m-1, 0)$-form $u_1 \wedge \cdots \wedge u_{m-1}$.

$\bigwedge^m W_{\alpha}$ splits as the direct sum of two irreducible
$\germg_{\alpha,\bC}$-submodules, $U_{1,\alpha}$ and $U_{2,\alpha}$,
having highest weights $2\mu_{m-1}$ and $2\mu_m$, respectively.
The space of $(m,0)$-forms in $\bigwedge^m W_{\alpha}$ is the
$1$-dimensional vector space $\bigwedge^m W_{\alpha}^{1,0}$ spanned by
$$u := u_1 \wedge \cdots \wedge u_m.$$
Since $u$ belongs to $U_{2,\alpha}$ (it is a vector of highest
weight), we conclude that $U_{2,\alpha}$ contains an $(m,0)$-form,
while $U_{1,\alpha}$ contains no $(m,0)$-forms.

Since $\disc T$ is not a square, the Tits index of $\widehat{G}$ is
$^2D_{m,r}^{(2)}$, where $r$ is the $k$-rank of $\widehat{G}$
\cite[p.~57]{Tits}.
In this case, $\Gal(\bar{k}/k)$ transposes the roots $\alpha_{m-1}$
and $\alpha_m$, while leaving $\alpha_1, \dots, \alpha_{m-2}$
invariant.
This implies that the weights $2\mu_{m-1}$ and $2\mu_m$ belong to the
same orbit of $\Gal(\bar{k}/k)$. Therefore, if $U$ is a $G$-module
such that $U_{\bC}$ contains $U_{1,\alpha}$, then $U_{\bC}$ must
contain a $G_\alpha (\bC)$-module equivalent to $U_{2,\alpha}$. It
follows that if $U$ is the smallest Hodge structure containing
$U_{1,\alpha}$, then  $U_{\bC}$ contains $(m,0)$-forms.
We may then proceed as in the proof of Case~3 of \cite[Theorem 5.1,
pp.~352--353]{Abdulali2} to show that any power of $A$ is dominated
by the set of all powers of $A$, and hence the usual Hodge conjecture
for all powers of $A$ implies the general Hodge conjecture for all
powers of $A$.
To complete the proof we note that by \cite[Theorem 4.1,
p.~673]{Abdulali3} the usual Hodge conjecture for $A$ implies the
usual Hodge conjecture for all powers of $A$.
\enddemo

{\it Remark} 4.2.
Murty \cite[Theorem~3.2, p.~203]{Murty} has shown that the Hodge ring
of an abelian variety of type III is never generated by divisors.
In \cite[Theorem~4.1, p.~673]{Abdulali3} we showed that if $A$ is a
sufficiently general abelian variety of type III, then the Hodge ring
of any power of $A$ is generated by divisors and Weil cycles.
Using Schoen's results \cite{Schoen} on the algebraicity of certain
Weil cycles, it is then possible to prove the usual Hodge conjecture
for all powers of certain type III abelian varieties of dimensions
$4$ \cite[Example 5.2, p.~674]{Abdulali3}, $6$ \cite[Example~5.1,
p.~674]{Abdulali3}, and, $8$ (van Geemen and Verra
\cite{GeemenVerra}).
In the $6$-dimensional case we were able to deduce the general Hodge
conjecture.
In the next corollary we do the same for the $4$-dimensional case.
The general Hodge conjecture in the $8$-dimensional case remains open.

\specialnumber{4.3}
\proclaim{{C}orollary}
\label{fourdim}
If $A$ is a simple $4$\/{\rm -}\/dimensional abelian variety of type {\rm III} whose
endomorphism algebra contains a square root of $-3${\rm ,} then{\rm ,} the
general Hodge conjecture is true for all powers of $A${\rm .}
\endproclaim

\demo{Proof}
By Moonen and Zarhin \cite[6.1, p.~568]{MoonenZarhin}, $A$ is a
general member of a PEL-family. By Shimura \cite[Proposition 17,
p.~180]{Shimura1}, the discriminant of the skewhermitian form
determined by a polarization of $A$ is not a square. As in
\cite[Example 5.2, p.~674]{Abdulali3}, the usual Hodge conjecture for
$A$ follows from \cite[Theorem 4.1, p.~673]{Abdulali3}, together with
Schoen's proof of the algebraicity of the required Weil cycles
\cite{Schoen}. Theorem \ref{nonsquare} now implies the general Hodge
conjecture for all powers of $A$.
\enddemo

\section{An extraordinary Hodge structure}

Throughout this section we retain the notations introduced in
previous sections, including the proof of Theorem \ref{nonsquare}.
First, we prove a converse to Theorem \ref{nonsquare}.

\proclaim{Proposition}
\label{squareprop}
Let $A$ be a simple abelian variety of type {\rm III.}
Assume that $A$ is a general member of a {\rm PEL-}\/family{\rm .}
Let $T$ be the skewhermitian form determined by a polarization of $A${\rm .}
If the discriminant of $T$ is a square{\rm ,} then $A$ is not dominated by
the set of powers of itself{\rm .}
\endproclaim

\demo{Proof}
The Tits index of $\widehat{G}$ is ${}^1D_{m,r}^{(2)}$
\cite[pp.~56--57]{Tits}, where $r$ is the $k$-rank of $\widehat{G}$.
In this case the $\star$-action (or $\Delta$-action) of
$\Gal(\overline{k}/k)$ on the simple roots of $\widehat{G}$ is
trivial.
By the arguments of Satake \cite[\S 1.1, pp.~218--220]{Satake2} this
implies that given any irreducible representation $\rho$ of
$\widehat{G}_{\overline k}$, some multiple $n\rho$ of $\rho$ is
defined over $k$.
Therefore there exists an irreducible representation
$\widehat{G} \to {\rm GL}(\widehat{U}_1)$
equivalent to a multiple of the representation with highest weight $2
\mu_{m-1}$.
Restricting scalars to~$\bQ$, we obtain an irreducible representation
$\chi$ of $G$ which is equivalent over $\bC$ to a multiple of
$\mathbold{\oplus}_{\alpha \in S} U_{1,\alpha}$.

Choose $\alpha \in S$, and let $M$ be the smallest Hodge substructure
of $\bigwedge^mV = H^m(A, \bQ)$ such that $M_{\bC}$ contains
$U_{1,\alpha}$.
Then the representation of $G$ on $M$ is equivalent to $\chi$.
We have seen in the proof of Theorem~\ref{nonsquare} that
$U_{1,\alpha}$ contains no $(m,0)$-forms, for any $\alpha \in S$.
Therefore $M_{\bC}$ contains no $(m,0)$-forms, and $M(1)$ is an
effective Hodge structure.
By Lemma \ref{notappear}, the representation of $G_\alpha (\bC)$ with
highest weight $2\mu_{m-1}$ does not occur in
$\bigwedge^{m-2}(W^n_\alpha)$ for any positive integer $n$, so $M(1)$
does not occur in the cohomology of any power of~$A$.
Thus $A$ is not dominated by the set of powers of itself.
\enddemo

\numbereddemo{{R}emark}
The general Hodge conjecture implies that the Hodge structure $M$
constructed in the above proof is supported on a divisor.
Since $\dim A = 2m[k:\bQ] > m$, the weak Lefschetz theorem implies
that $H^{m-2}(Y,\bQ) \cong H^{m-2}(A,\bQ)$ for a smooth hyperplane
section $Y$.
Since $M \subset H^m(A, \bQ)$, and we have seen that $H^{m-2}(A,\bQ)$
does not contain a Hodge substructure isomorphic to $M(1)$, we
conclude that $M$ is not supported on any smooth divisor.
\enddemo

\proclaim{Theorem}
\label{d4}
Let $A$ be an abelian variety satisfying the hypotheses of
Proposition {\rm \ref{squareprop}.} Assume further that $\dim_{D(A)} H_1(A,
\bQ) = 4${\rm .} Then there exists an abelian variety $B$ such that each
power of $A$ is dominated by the set of abelian varieties of the form
$A^i \times B^j${\rm .}
\endproclaim

\demo{Proof}
Let $\widetilde{G} \to G$ be the universal covering of $G$. Then
\eqref{gproduct} shows that $\widetilde{G}(\bR) \cong \prod_{\alpha
\in S} \widetilde{G}_{\alpha}$, where each $\widetilde{G}_{\alpha}$
is isomorphic to the real Lie group $\Spin(6,2)$. Since the
derivative of $\widetilde{G}(\bR) \to G(\bR)$ is an isomorphism, we
identify their Lie algebras.

Let $\sA \to \sV$ be a Kuga fiber variety obtained from the
$H_2$-morphism $G \to {\rm Sp}(V, \beta)$. This is the Hodge family
determined by $A$, so $A$ is the fiber $\sA_P$ for a general point $P
\in \sV$. Let $\rho$ be the composite $H_2$-morphism
$$\widetilde{G} \to G \to {\rm Sp}(V, \beta).$$
Satake (\cite[p.~451]{Satake1}, \cite{Satake2}) also constructed an
$H_2$-morphism
$$\rho' : \widetilde{G} \to {\rm Sp}(V', \beta')$$
such that the scalar extension of $\rho'$ to $\bR$ is equivalent to
$\mathbold{\oplus}_{\alpha \in S} \rho_{\alpha}$, where $\rho_{\alpha}$ is the
direct sum of two copies of the spin representation of
$\widetilde{G}_{\alpha}$ having highest weight $\mu_3$.
Let $\sB \to \sV$ be the Kuga fiber variety obtained from this
representation, and let $B := \sB_P$.
Proposition \ref{h2} implies that the Hodge group of $B$ is
$\rho'(\widetilde{G})$. The representation of
$\widetilde{G}_{\alpha}$ with highest weight $2\mu_3$ then appears in
$H^2(B^2, \bR)$.
Let $M'$ be an irreducible Hodge substructure of $H^2(B^2, \bQ)$
containing this representation.

The representation
$$\rho \oplus \rho' : \widetilde{G} \to {\rm Sp}(V \oplus V', \beta \oplus \beta')$$
also satisfies the $H_2$-condition.
The Kuga fiber variety obtained from this is the fiber product $\sA
\times_{\sV} \sB \to \sV$. The fiber over the point $P$ is $A \times
B$, so Proposition \ref{h2} implies that $G(A \times B) = (\rho
\oplus \rho')(\widetilde{G})$.
Since $M$ and $M'$ are equivalent $\widetilde{G}$-modules in the
cohomology of $A \times B^2$, with $M \subset H^4(A \times B^2, \bQ)$
and $M' \subset H^2(A \times B^2, \bQ)$, we have $M' \cong M(1)$.
The arguments in Case 3 of the proof of \cite[Theorem 5.1,
pp.~352--353]{Abdulali2} then show that any power of $A$ is dominated
by the set of abelian varieties of the form $A^i \times B^j$.
\enddemo

\numbereddemo{{R}emark}
The symmetric domains belonging to the hermitian groups
${\rm SO}^{\star}(8)$ and $\Spin(6,2)$ are isomorphic (cf. \cite[Exercise
1, p.~289]{Satakebook}).
This exceptional isomorphism is the reason why the case $m=4$ needs
separate treatment.
For odd $m$, $\disc T$ is not a square (it is negative), so Theorem
\ref{nonsquare} applies.
As explained in the proof of Corollary \ref{fourdim}, the case $m=2$
is also covered by Theorem \ref{nonsquare}.
The next theorem deals with the remaining cases.
\enddemo

\proclaim{Theorem}
\label{maintheorem}
Let $A$ be an abelian variety satisfying the hypotheses of
Proposition {\rm \ref{squareprop}.}
Assume further that $\dim_{D(A)} H_1(A, \bQ) > 4${\rm .}
Then $A$ is not dominated by the class of all abelian varieties{\rm .}
\endproclaim

\demo{Proof}
Lemma~\ref{quotient} shows that $G(M)$, the Hodge group of the Hodge
structure $M$ constructed in the proof of Proposition
\ref{squareprop}, is a quotient of~$G$. Since $G$ is (almost) simple,
and $M$ is a nontrivial Hodge structure, we see that $G(M)$ is
isogenous to~$G$. Suppose $M(1)$ occurs as a Hodge substructure of
$H^{m-2}(B, \bQ)$ for an abelian variety~$B$. Clearly, the Hodge
groups of $M$ and $M(1)$ are equal. Using Lemma~\ref{quotient} again,
we find that $G(M)$ is a quotient of $G(B)$. Therefore $\germg (B) =
\germg \oplus \germg_2$, where $\germg_2$ is a reductive Lie algebra
which acts trivially on $M(1)$. Since $\germg$ is simple, it follows
that $\germg^{\der}(B) = \germg \oplus \germg_2^{\der}$.

Let $W := H_1(B, \bQ)$, and let $\gamma$ be a Riemann form for $B$.
As explained in \S 2.3, the inclusion $\rho :
G^{\der}(B) \to {\rm Sp}(W, \gamma)$ is an $H_1$-morphism which defines a
Hodge family $\sB \to \sW$ having $B$ as a fiber over a general point
$Q \in \sW$.
We have the primary decomposition $W = \mathbold{\oplus}_{\alpha} W_{\alpha}$,
where the $W_\alpha$ are the maximal primary $G^{\der}(B)$-submodules
of $W$ (\S 2.4). Since $G(\bR)$ has no compact factors,
Lemma~\ref{primarylemma} implies that any $W_{\alpha}$ is either
trivial as a $\germg$-module or trivial as a $\germg_2^{\der}$-module.

Let $W_1$ be the sum of those $W_{\alpha}$'s on which $\germg$ acts
nontrivially, and let $W_2$ be the sum of those $W_{\alpha}$'s on
which $\germg$ acts trivially.
Then $W = W_1 \oplus W_2$, where $\germg^{\der}_2$ acts trivially
on~$W_1$, and $\germg$ acts trivially on $W_2$.
The $H_1$-morphisms
$$\rho_i : G^{\der}(B) \to {\rm Sp}(W_i, \gamma | W_i \times W_i), \qquad i = 1,2,$$
define Kuga fiber varieties $\sB_i \to \sW$ whose fiber product is
isogenous to $\sB \to \sW$.
We then have $B$ isogenous to $B_1 \times B_2$, where $B_i$ is the
fiber of $\sB_i$ over the point $Q$.
Since $\sB \to \sW$ is a Hodge family, a theorem of Mumford
\cite[p.~348]{Mumford2} implies that it contains a fiber of CM-type.
It follows that $\sB_i \to \sW$ also contain CM fibers, and the same
theorem of Mumford implies that they are Hodge families.
Since $Q$ is a general point of $\sW$, we have $G^{\der}(B_i) =
\rho_i (G^{\der}(B))$. Thus $\germg(B_1) = \germg$, and
$\germg^{\der}(B_2) = \germg^{\der}_2$.
Since $\germg_2$ acts trivially on $M(1)$, we see that $M(1)$ occurs
in the cohomology of $B_1$, and we assume without loss of generality
that $B = B_1$, and $G(B)$ is isogenous to $G$.

The isomorphism of $G$-modules
$$ H^m(A, \bQ) \supset M \stackrel{\cong}{\longrightarrow} M(1)
\subset H^{m-2}(B, \bQ)$$
is induced by a Hodge cycle $\zeta$ on $A \times B$.
Now, $G(A \times B)$ is a subgroup of $G(A) \times G(B)$ whose
projection to each factor is surjective.
By Goursat's Lemma (cf.~Gordon \cite[Proposition~B.71,
p.~316]{Gordon}), $\germg (A \times B)$ is either $\germg(A) \times
\germg(B)$, or the graph of an isomorphism from $\germg(A)$ to
$\germg(B)$.
Since $\zeta$ is invariant under $G(A \times B)$, the possibility
$\germg(A) \times \germg(B)$ is ruled out, so we conclude that $G(A
\times B)$ is isogenous to $G$.

Satake's classification of $H_1$-morphisms \cite{Satake1} shows that
the inclusion of $G(A \times B)$ into the symplectic group of $H_1(A
\times B, \bQ)$ satisfies the $H_2$-condition and is a multiple of
the standard representation.
Lemma \ref{rigidity} now implies that $B$ is isogenous to a power of
$A$, contradicting Proposition \ref{squareprop}.
\enddemo
 
 \end{document}